\documentclass[11pt]{amsart}
\usepackage{geometry} 
\usepackage{url}               % See geometry.pdf to learn the layout options. There are lots.
\usepackage{graphicx}
\usepackage{amssymb}
\usepackage{hyperref}
\usepackage[shortalphabetic]{amsrefs}
\usepackage{mathptmx}   
\newtheorem{theorem}{Theorem}
\newtheorem{proposition}[theorem]{Proposition}

{\theoremstyle{remark}

}

\newcommand{\RR}{\mathbb R}

\DeclareMathOperator{\tr}{tr}

\title[Contracting surfaces in space]{Moving surfaces by non-concave
curvature functions}
\author{Ben Andrews}
\address{MSI, ANU, ACT 0200 Australia}
\email{Ben.Andrews@anu.edu.au}
\thanks{Research supported by Discovery grants DP0344221
and DP0985802 of the Australian Research Council.}
\begin{document}
\begin{abstract}
A convex surface contracting by a strictly monotone,
homogeneous degree one function of its principal curvatures remains smooth 
until it contracts to a point in finite time, and is asymptotically 
spherical in shape.  No assumptions are made on the concavity of the
speed as a function of principal curvatures.  We also discuss motion by
functions homogeneous of degree greater than $1$ in the principal curvatures.
\end{abstract}

\maketitle

\section{Introduction}\label{sec:intro}

% history of homogeneous degree one contraction

Several authors have considered
convex hypersurfaces contracting by homogeneous degree one symmetric
functions of their principal curvatures:  Huisken
\cite{Hu} proved that the mean curvature flow
contracts such hypersurfaces to a point in finite time while
making their shape spherical.  Chow proved similar
results for $n$-dimensional hypersurfaces moving by the $n$-th root
of Gauss curvature \cite{Ch1} and (with additional convexity
assumptions on the initial data) by the square root of
scalar curvature \cite{Ch2}.  The author treated a large family of 
such equations \cites{A1,AndrewsPinching}, satisfying some requirements
of concavity of the speed in the principal curvatures.

Since the
result in \cite{A1} holds for both  concave and convex functions of the
principal curvatures, it seems possible that these convexity
assumptions could be significantly weakened.  This paper confirms that no concavity assumptions are needed for 
surfaces moving in space:

% the result

\begin{theorem}\label{thm:result}
	Let $F$ be a smooth, symmetric, homogeneous degree 1 function
	$F$ defined on the positive quadrant in $\RR^2$, with strictly
	positive derivative in each argument, normalized to have $F(1,1)=1$,  Then 
	for any smooth, strictly convex surface $M_0=x_0(S^2)\subset\RR^3$ there is a
	unique family of smooth, strictly convex surfaces 
	$\{M_t=x_t(S^2)\}_{0\leq t<T}$ satisfying
	\begin{equation}\label{eq:flow}
		\begin{split}
			{\partial x\over\partial t}(z,t)&=-F(\kappa_1(z,t),\kappa_2(z,t))
			\nu(z,t)\\
			x(z,0)&=x_0(z)\\
		\end{split}
	\end{equation}
	where $\nu(z,t)$ is the outward normal and
	$\kappa_1(z,t)$ and $\kappa_2(z,t)$ are the principal curvatures
	of $M_t$.  $M_t$ converges uniformly to a point
	$p\in\RR^3$ as $t$ approaches $T$, while the rescaled maps
	${x_t-p\over\sqrt{2(T-t)}}$ converge smoothly to an
	embedding $\tilde x_T$ with image equal to the unit sphere about the
	origin.
\end{theorem}
% special aspects of the two-dimensional case

The proof of this result follows the same basic framework as the papers
mentioned above:  The crucial step is to obtain a bound on the
`pinching ratio', which is the supremum over the surface of the
ratio of largest to smallest principal curvatures at each point. This is
obtained using a maximum principle
argument. In the previous papers the concavity assumptions on the speed
came into the argument at two points:  First, to obtain terms of a
favourable sign in the argument to control the pinching ratio, and second
to allow the application of the second derivative H\"older estimates
of Krylov \cite{K}.  In the two-dimensional case the author recently
proved second derivative H\"older estimates which apply without any
assumption of concavity \cite{A2}.  The main new estimate of this paper
gives a bound on the pinching ratio without requiring any
concavity of the speed function, thus removing all concavity requirements
from the proof.

The last section of the paper discusses flows in which the
speed function is homogeneous of degree greater than 1 in the 
principal curvatures:  In particular, it is shown that any parabolic flow
with speed homogeneous of degree $\alpha>1$ preserves pinching ratios
which are less than or equal to a critical value $r_0(\alpha)$ which 
depends only on $\alpha$.  On the other hand there are surfaces with
pinching ratio as close to $r_0(\alpha)$ as desired, for which the
pinching ratio becomes worse under the flow.  

% expansion flows?

\section{Notation and preliminary results}\label{sec:prelim}

Suppose that the initial surface $M_0$ is given by
a smooth embedding $x_0: S^2\to\RR^3$.  The aim is to construct a smooth
family of embeddings $x: S^2\times[0,T)\to\RR^3$ satisfying the 
evolution equation.

The existence of a smooth solution for a short time is guaranteed since
the flow is equivalent to a scalar, strictly parabolic equation (for
example one can write the evolving surfaces as graphs over a sphere).

By convention the unit normal $\nu$ points in the outward
direction for a convex surface $M$, and the principal curvatures
are positive.  Choose local coordinates $y^1,y^2$ about any point
$z$ such that the tangent vectors $e_i={\partial x\over\partial y^i}$
are orthonormal at $z$.  Then the principal curvatures are the
eigenvalues of the second fundamental form, which is the symmetric
bilinear form defined by 
\[
	h_{ij}=-\left\langle {\partial^2x\over\partial y^i\partial
	y^j},\nu\right\rangle = \left\langle {\partial\nu\over\partial
	y^i},{\partial x\over \partial y^j}\right\rangle.
\]

The covariant derivatives of a tangent vector field $X$ on $M$ in
the direction of the tangent vector $e_j$ is given by the expression 
\[
	\nabla_{e_j}X=\pi\left({\partial X\over\partial y^j}\right),
\]
where $\pi$ is the orthogonal projection on the tangent space of $M$. The
covariant derivative of the second fundamental form is the tensor defined
by
\[
	\nabla_ih_{jk} = {\partial h_{jk}\over\partial y^i}
	-h\left(e_j,\nabla_{e_i}e_k\right)
	-h\left(\nabla_{e_i}e_j,e_k\right).
\]
The Codazzi identity says that this is totally symmetric.

The function $F$ is assumed to be a smooth symmetric function of the
principal curvatures, and so can also be written as a smooth function
of the elementary symmetric functions of the principal curvatures  \cite{Glaeser}, hence also of 
the components of the second fundamental form.  We remark that a 
Theorem of Schwarz \cite{Schwarz} guarantees that conversely, any smooth $SO(n)$-invariant
function of the components of the second fundamental form can be written as a
smooth symmetric function of the principal curvatures.  Define
$\dot F^{ij}={\partial F\over\partial h_{ij}}$ and $\ddot F^{klmn}=
{\partial^2 F\over\partial h_{kl}\partial h_{mn}}$.  The monotonicity of
$F$ as a function of each principal curvature implies that $\dot F$ is
a positive definite symmetric matrix.

Suppose $G$ is any other symmetric function of the principal curvatures,
and define $\dot G^{ij}$ and $\ddot G^{klmn}$ to be the first and second
derivatives of $G$ with respect to the components of the second
fundamental form.  It was shown in \cite{A1} that if $F$ is homogeneous
of degree 1 then $G$ evolves according to the following evolution
equation:
\begin{equation}\label{eq:curvevol}
	\begin{split}
		{\partial G\over\partial t}
		&=\dot F^{ij}\nabla_i\nabla_jG
		+(\dot G^{ij}\ddot F^{klmn}-\dot F^{ij}\ddot G^{klmn})\nabla_ih_{kl}
		\nabla_jh_{mn}\\
		&\quad\null+\dot G^{ij}h_{ij}\dot
		F^{kl}g^{mn}h_{km}h_{ln}
	\end{split}
\end{equation}
where the sum is over repeated indices.  In particular,
\begin{equation}\label{eq:speedevol}
	{\partial F\over\partial t} = \dot F^{ij}\nabla_i\nabla_jF+F\dot
	F^{ij}g^{kl}h_{ik}h_{jl}.
\end{equation}
The parabolic maximum principle therefore implies that the minimum
of $F$ over the surface $M_t$ is non-decreasing in time as long as the
solution remains smooth and convex.

\section{Estimate on the pinching ratio}\label{sec:pinch}

In this section the main estimate of the paper is proved by applying the
maximum principle to the evolution equation for the quantity 
\begin{equation}\label{eq:defG}
	G={(\kappa_2-\kappa_1)^2\over (\kappa_2+\kappa_1)^2}=
	{2|h|^2-(\tr h)^2\over (\tr h)^2}.
\end{equation}
Since $G$ is homogeneous of degree zero, the Euler relation gives
$\dot G^{ij}h_{ij}=0$, and the last term in equation \eqref{eq:curvevol} 
vanishes.

Consider a point $(z,t)$ where $G$ attains a spatial maximum at some time
$t$ in the interval of existence of the smooth convex solution of the
flow. Consider the second term in Equation
\eqref{eq:curvevol} at such a point.  This is a quadratic term in the 
components of $\nabla h$.  The Codazzi identity implies that there are 
four distinct components of
$\nabla h$,  so this term is defined by a $4\times 4$ matrix.  However,
at the maximum point the first derivatives of $G$
vanish, so two components of $\nabla h$ can be eliminated,
leaving only a
$2\times 2$ matrix.  It will turn out that in an orthonormal basis diagonalising 
the second fundamental form at some point and time, this matrix is also
diagonal and can be easily understood.

The computation of these terms will require results from \cite{A1} which
describe the components of $\ddot F$ and $\ddot G$ in an orthonormal
frame in which $h_{ij}=\text{\rm diag}(\kappa_1,\kappa_2)$:
\begin{equation*}
	\begin{split}
		\ddot F^{11,11}&={\partial^2F\over\partial\kappa_1^2}\\
		\ddot F^{11,22}=\ddot
		F^{22,11}&={\partial^2F\over\partial\kappa_1\partial\kappa_2}\\
		\ddot F^{22,22}&={\partial^2F\over\partial\kappa_2^2}\\
		\ddot F^{12,12}=\ddot
		F^{21,21}&={{\partial F\over\partial\kappa_1}-{\partial
		F\over\partial\kappa_2}\over\kappa_1-\kappa_2}\\
	\end{split}
\end{equation*}
The last of these identities is to be interpreted as a limit if
$\kappa_1=\kappa_2$.  It follows that the terms in the evolution
equation for $G$ are as follows in a frame diagonalising the second
fundamental form:
\begin{equation*}
	\begin{split}
		Q&=
		\left(\dot G^{ij}\ddot F^{klmn}-\dot F^{ij}\ddot G^{klmn}\right)
		\nabla_i h_{kl}\nabla_jh_{mn}\\
		&= \left({\partial G\over\partial\kappa_1}
		{\partial^2F\over\partial\kappa_1^2}
		-{\partial F\over\partial\kappa_1}
		{\partial^2G\over\partial\kappa_1^2}\right)
		(\nabla_1h_{11})^2
		+\left({\partial G\over\partial\kappa_1}
		{\partial^2F\over\partial\kappa_2^2}
		-{\partial F\over\partial\kappa_1}
		{\partial^2G\over\partial\kappa_2^2}\right)(\nabla_1h_{22})^2\\
		&\quad\null
		+2\left({\partial G\over\partial\kappa_1}
		{\partial^2F\over\partial\kappa_1\partial\kappa_2}
		-{\partial F\over\partial\kappa_1}
		{\partial^2 G\over\partial\kappa_1\partial\kappa_2}\right)
		\nabla_1h_{11}\nabla_1h_{22}\\
		&\quad\null+\left({\partial G\over\partial\kappa_2}
		{\partial^2F\over\partial\kappa_1^2}
		-{\partial F\over\partial\kappa_2}
		{\partial^2G\over\partial\kappa_1^2}\right)
		(\nabla_2h_{11})^2
		+\left({\partial G\over\partial\kappa_2}
		{\partial^2F\over\partial\kappa_2^2}
		-{\partial F\over\partial\kappa_2}
		{\partial^2G\over\partial\kappa_2^2}\right)(\nabla_2h_{22})^2\\
		&\quad\null
		+2\left({\partial G\over\partial\kappa_2}
		{\partial^2F\over\partial\kappa_1\partial\kappa_2}
		-{\partial F\over\partial\kappa_2}
		{\partial^2 G\over\partial\kappa_1\partial\kappa_2}\right)
		\nabla_2h_{11}\nabla_2h_{22}\\
		&\quad\null+2
		{{\partial G\over\partial\kappa_1}
		{\partial F\over\partial\kappa_2}-
		{\partial G\over\partial\kappa_2}
		{\partial F\over\partial\kappa_1}
		\over \kappa_2-\kappa_1}(\nabla_1h_{12})^2
		+2{{\partial
		G\over\partial\kappa_1} {\partial F\over\partial\kappa_2}-
		{\partial G\over\partial\kappa_2}
		{\partial F\over\partial\kappa_1}
		\over \kappa_2-\kappa_1}(\nabla_2h_{12})^2
	\end{split}
\end{equation*}
At a maximum point of $G$, 
$G$ is non-zero (otherwise $M_t$ is a sphere and the proof is trivial)
and it be can assumed without loss of generality that $\kappa_2>\kappa_1$.
The gradient conditions on $G$ then give two equations (since
$\dot G$ does not vanish except when $G=0$):
\begin{equation*}
		\nabla_1h_{11}=-{{\partial G\over\partial\kappa_2}
		\over{\partial G\over\partial\kappa_1}}\nabla_1h_{22};\qquad
		\nabla_2h_{22}=-{{\partial G\over\partial\kappa_1}
		\over{\partial G\over\partial\kappa_2}}\nabla_2h_{11}.
\end{equation*}
The degree-zero homogeneity of $G$ implies by the Euler relation that
$\kappa_1{\partial G\over\partial\kappa_1}+
\kappa_2{\partial G\over\partial\kappa_2}=0$.  Homogeneity also implies the identity
\[
	\kappa_1^2{\partial^2G\over\partial\kappa_1^2}
	+\kappa_2^2{\partial^2G\over\partial\kappa_2^2}
	+2\kappa_1\kappa_2{\partial^2 G\over\partial\kappa_1\partial\kappa_2}=0.
\]
Similarly, the degree 1 homogeneity of $F$ gives the following identities:
\[
		{\partial^2F\over\partial\kappa_1^2}=-{\kappa_2\over\kappa_1}
		{\partial^2F\over\partial\kappa_1\partial\kappa_2};\qquad
		{\partial^2F\over\partial\kappa_2^2}=-{\kappa_1\over\kappa_2}
		{\partial^2F\over\partial\kappa_1\partial\kappa_2};\qquad
		\kappa_1{\partial F\over\partial\kappa_1}
		+\kappa_2{\partial F\over\partial\kappa_2}=F.
\]
Substituting these expression into the expression for $Q$ above and
applying the Codazzi symmetries $\nabla_1h_{12}=\nabla_2h_{11}$ and
$\nabla_2h_{12}=\nabla_1h_{22}$, one finds that all of the terms
involving second derivatives of $F$ and $G$ disappear, leaving
\[
	Q={2F{\partial G\over\partial\kappa_1}\over \kappa_2(\kappa_2-\kappa_1)}
	(\nabla_1h_{22})^2+
	{2F{\partial G\over\partial\kappa_1}\over \kappa_2(\kappa_2-\kappa_1)}
	(\nabla_2h_{11})^2.
\]
Now observe that ${\partial G\over\partial\kappa_1}=
-{4\kappa_2(\kappa_2-\kappa_1)\over(\kappa_1+\kappa_2)^3}<0$, 
and therefore $Q\leq 0$ and ${\partial G\over\partial t}\leq 0$ at the
maximum point.  Therefore the supremum $\bar G$ of $G$ is non-increasing 
in time,  and the pinching ratio $r={2\over1-\sqrt{\bar G}}-1$ is also 
non-increasing in time.
This proves the following:

\begin{proposition}\label{prop:maxpr}
	If $x: S^2\times [0,T)\to{\mathbb R}^3$ is a smooth family of convex embeddings
	satisfying  \thetag{1}, then for $t\in(0,T)$ the pinching
	ratio of
	$M_t=x_t(S^2)$ is no greater than the pinching ratio of $M_0$.
\end{proposition}

If follows from \cite{A1}*{Lemma 5.4} that the surfaces
$M_t$ have bounded ratio of circumradius $r_+$ to inradius $r_-$:  There exists $C_0$ such that $r_+(M_t)\leq C_0r_-(M_t)$ for all $t$ in the interval of existence.

\section{Regularity and convergence}\label{sec:reg}

Now we discuss the proof of convergence to a point and the limit under rescaling.  The key step is
to derive estimates on curvature and its higher derivatives.  To achieve the correct dependence of the estimates on the geometry we rescale to bring the hypersurfaces to a fixed size:

For each $t_0\in[0,T)$, let $p_{t_0}$ be an incentre of $M_{t_0}$, and let 
$x_{t_0}$ be the solution of Equation \eqref{eq:flow} defined by
$x_{t_0}(z,t) = r_-(M_{t_0})^{-1}\left(x(z,t_0+r_-(M_{t_0})^2t)-p_{t_0}\right)$ for $t\in [-t_0/r_-^2,(T-t_0)/r_-^2)$.  Denote the hypersurface
$x_{t_0}(S^2,t)$ by $M_{t_0,t}$.  Then $M_{t_0,0}$ lies outside the unit ball about the origin in $\RR^3$, and inside the ball of radius $C_0+1$, for every $t_0$.  By the comparison principle, 
$M_{t_0,t}$ lies outside the ball of radius $\sqrt{1-2t}$ and inside the ball of radius $\sqrt{(C_0+1)^2-2t}$, for each $t\in[0,1/2)$ in the interval of existence.

The pinching estimate of Proposition \ref{prop:maxpr} implies that the eigenvalues of $\dot F$ are
bounded above and below by positive constants:  By homogeneity we have 
$\frac{\partial F}{\partial\kappa_i}(\kappa_1,\kappa_2) =\frac{\partial F}{\partial\kappa_i}\left(\frac{\kappa_1}{\kappa_1+\kappa_2},\frac{\kappa_2}{\kappa_1+\kappa_2}\right)$, so the supremum and infimum are attained on the compact set $\left\{(a,1-a):\ |a-\frac12|\leq \frac{\sqrt{C}}{2}\right\}$, and hence are finite and positive respectively.

Bounds above on $F$ on the rescaled hypersurfaces follow from the argument of Tso \cite{Tso} as presented in \cite{A1}*{Theorem 7.5}.  Together with the pinching estimate this implies a uniform upper bound on the principal curvatures of $M_{t_0,t}$ for any $t_0$ with $0\leq t\leq \frac14$.

For any orthonormal basis for $\RR^3$, and each $t_0\in[0,T)$ and $t\in[0,\frac14]$ with $t_0+r_-(t_0)^2t<T$,
$$
M_{t,t_0}\cap\left\{(x,y,z)\in\RR^3: \ z<0,\ x^2+y^2<\frac{1}{16}\right\}
= \left\{(x,y,u_{t_0}(x,y,t):\ x^2+y^2<\frac{1}{16}\right\},
$$
where $-(C_0+1)<u_{t_0}(x,y,t)<0$ and 
$|Du_{t_0}|(x,y,t)\leq 4(C_0+1)$, and
\begin{equation}\label{eq:grapheqn}
\frac{\partial u_{t_0}}{\partial t} = \tilde F(D^2u_{t_0},Du_{t_0})=F\left(\left(g(Du_{t_0})\right)^{-1/2}(D^2u_{t_0})\left(g(Du_{t_0})\right)^{-1/2}\right),
\end{equation}
where 
$g(V) = I+V^TV$.  The matrix $g$ has eigenvalues bounded below by $1$ and above by $1+16(C_0+1)^2$, and so the derivatives of $\tilde F$ with respect to the components of $D^2u$ are comparable to the derivatives of $F$, which are given by $\frac{\partial F}{\partial\kappa_i}$, and Equation \eqref{eq:grapheqn} is uniformly parabolic.  The Krylov-Safonov Harnack estimate \cite{K} gives positive lower bounds for $F$ (hence also positive lower bounds for all principal curvatures) for $M_{t,t_0}$.  Theorem 5 of \cite{A2} gives uniform H\"older bounds on the second spatial derivatives of $u_{t_0}$ on the same range.  Finally, Schauder estimates \cite{Lieb}*{Theorem 4.9} give uniform bounds on all higher derivatives of $u_{t_0}$.  Since the choice of basis was arbitrary, these imply uniform bounds on curvature and its higher derivatives on $M_{t_0,t}$.  

It follows that the hypersurfaces $M_{t_0,t}$ extend to exist on $S^2\times [0,\frac14)$ for every $t_0\in[0,T)$, and consequently $T>t_0+r_-(M_{t_0})^2/4$ for all $t_0<T$, so that $r_-(M_{t_0})^2\leq 4(T-t_0)$ for all $t_0<T$, so the inradius $r_-(M_t)$ (hence also the circumradius $r_+(M_t)$) approach zero as $t$ approaches $T$.

Finally, the uniform bounds on the families of hypersurfaces $\{M_{t_0,t}\}$ imply that there exists a subsequence $t_k\to T$ such that the families $M_{t_k,t}$ converge smoothly to a limiting family $\bar M_t$, $0\leq t\leq \frac14$, which is again a solution of \eqref{eq:flow}, on which the pinching ratio is constant in time.  By the strong maximum principle applied to the evolution equation for $G$ derived in Section \ref{sec:pinch}, the limit solution is a shrinking sphere.  This proves sub-sequential convergence of the rescaled solutions to a sphere, and stronger convergence can be deduced by considering the linearization of the flow about the shrinking sphere solution as in \cite{ACLF}*{Propositions 40--41}.

\section{Remarks on higher degrees of homogeneity}\label{sec:highdeg}

In this section the methods of the previous sections are applied
to flows in which the speed is homogeneous of
some degree $\alpha>1$ in the principal curvatures.  The conclusion
is that such flows do not preserve large values of the pinching ratio.  In work of the author
on motion of surfaces by Gauss curvature \cite{A3}, it was
shown that the maximum difference between the principal curvatures does
not get any larger under this flow.  Thus the fact that the pinching
ratio does not improve does not rule out the possibility that other
curvature estimates may yield useful results.

The evolution of a curvature function $G$ under a flow with speed
$F$ which is homogeneous of degree $\alpha$ is as follows (compare
equation \eqref{eq:curvevol}):
\begin{equation}\label{eq:curvevol2}
	\begin{split}
		{\partial G\over\partial t}
		&=\dot F^{ij}\nabla_i\nabla_jG
		+(\dot G^{ij}\ddot F^{klmn}-\dot F^{ij}\ddot G^{klmn})\nabla_ih_{kl}
		\nabla_jh_{mn}\\
		&\quad\null+\dot G^{ij}h_{ij}\dot
		F^{kl}g^{mn}h_{km}h_{ln}-(\alpha-1)F\dot G^{ij}g^{kl}h_{ik}h_{jl}
	\end{split}
\end{equation}
If $G$ is homogeneous of degree zero, then the first term on the
second line vanishes.  Also, $\dot G^{ij}g^{kl}h_{ik}h_{jl}=
{\partial G\over\partial\kappa_1}\kappa_1^2+
{\partial G\over\partial\kappa_2}\kappa_2^2=
{\partial G\over\partial\kappa_2}\kappa_2(\kappa_2-\kappa_1)\geq0$
for $G$ as in Section \ref{sec:pinch} above.  Therefore the last term is 
non-positive.  It remains to understand the gradient terms, as before.

The difference from the computation in Section \ref{sec:pinch} arises 
from a change in the Euler identities:  
\[
		\kappa_1^2{\partial^2F\over\partial\kappa_1^2}
		+2\kappa_1\kappa_2{\partial^2F\over\partial\kappa_1\partial\kappa_2}
		+\kappa_2^2{\partial^2F\over\partial\kappa_2^2}=\alpha(\alpha-1)F;\qquad
		\kappa_1{\partial F\over\partial\kappa_1}
		+\kappa_2{\partial F\over\partial\kappa_2}=\alpha F.
\]
These lead to the
following expression for the gradient terms:
\[
	\begin{split}
		Q&=
		\left(\dot G^{ij}\ddot F^{klmn}-\dot F^{ij}\ddot G^{klmn}\right)
		\nabla_i h_{kl}\nabla_jh_{mn}\\ 
		&={\partial G\over\partial\kappa_1}
		\left(\alpha(\alpha-1){F\over\kappa_2^2} 
		+2\alpha {F\over\kappa_2(\kappa_2-\kappa_1)}
		\right)(\nabla_1h_{22})^2\\
		&\quad\null+{\partial G\over\partial\kappa_2}
		\left(\alpha(\alpha-1){F\over\kappa_1^2}
		-2\alpha{F\over\kappa_1(\kappa_2-\kappa_1)}
		\right)(\nabla_2h_{11})^2.
	\end{split}
\]
For this to be negative the pinching ratio $r=\kappa_2/\kappa_1$
must satisfy
 $2r+(\alpha-1)(r-1)\geq 0$ and 
$(\alpha-1)r-(\alpha-1)-2\leq 0$.  The first  is always
true since $r\geq 1$, but the second holds only if
\[
	r\leq r_0(\alpha)=1+\frac{2}{\alpha-1}.
\]
Thus the flow will improve pinching ratios no greater than $r_0(\alpha)$.  
The argument to prove smooth convergence to a sphere then follows exactly that in \cite{AM}.
The precise result is as follows:

\begin{theorem}\label{thm:highdeg}
Let $F$ be a smooth function defined on the positive cone, homogeneous
of degree $\alpha>1$, strictly increasing in each argument, and normalized to have $F(1,1)=1$.  Then
for any surface $M_0=x_0(S^2)$ which is smooth and strictly convex
with pinching ratio $r\leq r_0(\alpha)$ there exists a unique smooth
solution $x: S^2\times [0,T)\to{\mathbb R}^3$ of the evolution equation
\[
	\begin{split}
		{\partial x\over\partial
		t}(z,t)&=-F(\kappa_1(z,t),\kappa_2(z,t))\nu(z,t);\\
		x(z,0)&=x_0(z).
	\end{split}
\]
The surfaces $M_t=x_t(S^2)$ converge
to a point $p\in{\mathbb R}^3$ as $t\to T$, and the rescaled hypersurfaces
${M_t-p\over ((1+\alpha)(T-t))^{1/(1+\alpha)}}$ converge in $C^\infty$ to the unit sphere about the origin.
\end{theorem}

Next we show that this result cannot be improved, by constructing examples of smooth, strictly convex surfaces for
which the pinching ratio becomes larger, for any flow with a speed 
homogeneous of degree $\alpha>1$.  Remarkably, the particular surface which provides a counterexample
depends only on the degree of homogeneity $\alpha$.

Consider surfaces given by rotating the graph $y=u(x)$ about the $x$ axis.
The principal curvatures of this surface are given in terms of $u$ by
\[
	\kappa_1={1\over u\sqrt{1+(u')^2}}
\]
in the direction within the $y-z$
plane, and
\[
	\kappa_2=-{u''\over(1+(u')^2)^{3/2}}
\]
in the direction along the $x$ axis.  Therefore the
ratio of principal curvatures is equal to
\[
	r=-{uu''\over 1+(u')^2}.
\]

The ratio of principal curvatures can be prescribed as a function
of distance from the axis:   
Let $f:\RR\to\RR$ be any smooth positive even function with $f(0)=1$.
Then it is necessary to solve the following ordinary differential equation 
with initial data of the form $u(0)=U$, $u'(0)=1$:
\[
	f(u)=-{uu''\over(1+(u')^2}.
\]
This can be integrated to give
\[
	(u')^2=\exp\left\{2\int_u^U{f(z)\over z}\,dz\right\}-1,
\] 
which can in turn be integrated to give a smooth solution on $[0,L)$ with $u(x)\sim C\sqrt{L^2-x^2}$ as
$x$ approaches $L$.  Extending this to be even in $x$ gives
a smooth, strictly convex surface with ratio of principal curvatures
equal to $f(u)$ where $u$ is the distance from the axis of rotation.

Now choose $f(u)$ to be a smooth function with $f(u)=r_1>1$ for
$u\geq u_0$ and $f(u)\in(1,r_1)$ for $0<u<u_0$.  If $U>u_0$,
then this defines a smooth,  strictly convex surface with pinching ratio
equal to $r_1$, with this pinching ratio attained on an open annular
region away from the `poles' of the surface.  In this region,
$\nabla G=0$ and $\nabla\nabla G=0$, where $\nabla$ is the 
covariant derivative and $G$ is as given in section 3.  A direct
calculation also shows that $\nabla_1h_{22}=0$ everywhere,
while 
\[
	|\nabla_2h_{11}|^2= {(u')^2(1+(u')^2)\over u^4}(f(u)-1)^2
\]
which is certainly not identically zero on this annular region.  Therefore
in the evolution equation \eqref{eq:curvevol2} for $G$, the first term is zero,
the second term is positive, and the last term is negative.  To get further
one can write out these terms explicitly in terms of $u$:  At points with
$r=r_1$,
\[
	{\partial G\over\partial t}
	= {F{\partial G\over\partial\kappa_2}(r_1-1)\over r_1u^2(1+(u')^2)}
	\left(\alpha\left((\alpha-1)r_1(r_1-1)-2\right)(u')^2(1+(u')^2)^3
	-(\alpha-1)r_1^2\right)
\]

But now by choosing $u_0$ sufficiently close to zero while keeping $U$
fixed, it can be guaranteed that there are points achieving the pinching ratio
which have $u'$ as large as desired.  As long as 
$r_1>r_0(\alpha)$, the positive
first term in the bracket dominates the negative second term, so
that $G$ is strictly increasing, and the pinching ratio becomes
larger for small positive times.

\end{document}